\documentclass[10pt]{amsart}

\usepackage{latexsym}
\usepackage{labelfig}
\usepackage{amssymb}
\usepackage[cp850]{inputenc}
\usepackage{epsfig}
\usepackage{epstopdf}
\usepackage{psfrag}
\usepackage{amsthm}
\usepackage{amscd}
\usepackage{amsmath}
\usepackage{amsfonts}
\usepackage{graphics}

\newtheorem{thm}{Theorem}				\newtheorem{lem}{Lemma}
\newtheorem{cor}[thm]{Corollary}			\newtheorem{prop}[thm]{Proposition}
				
\newtheorem*{defi*}{Definition}	\newtheorem*{ex*}{Example}
\newtheorem*{thm*}{Theorem}			\newtheorem*{cor*}{Corollary}
\newtheorem*{rmk*}{Remark}			
		\newtheorem*{prop*}{Proposition}

\renewcommand{\theequation}{\thesection.\arabic{equation}}
\let\ssection=\section
\renewcommand{\section}{\setcounter{equation}{0}\ssection}

\newtheorem*{namedtheorem}{\theoremname}
\newcommand{\theoremname}{testing}
\newenvironment{named}[1]{\renewcommand{\theoremname}{#1}\begin{namedtheorem}}{\end{namedtheorem}}

\theoremstyle{remark}

		\newcommand{\CX}{\mathcal X}

\newcommand{\lra}{\longrightarrow}

\newcommand{\bs}{\backslash}		
\DeclareMathOperator{\arcsinh}{arcsinh}	
	
\DeclareMathOperator{\ii}{i}	

\begin{document}

\title[]{Small filling sets of curves on a surface}

\author[J.~W.~Anderson]{James W. Anderson}
\address{School of Mathematics \\
  University of Southampton \\
  Southampton, UK}
\email{j.w.anderson@soton.ac.uk}

\author[H.~Parlier]{Hugo Parlier}
\address{Department of Mathematics, University of Fribourg\\
Fribourg, Switzerland}
\email{hugo.parlier@gmail.com}

\author[A.~Pettet]{Alexandra Pettet}
\address{Department of Mathematics \\
  University of Michigan\\
  Ann Arbor, USA \newline \indent
Mathematical Institute \\ University of Oxford \\ Oxford, UK}
\email{alexandra.pettet@maths.ox.ac.uk}

\thanks{\tiny The first and third authors were supported by EPSRC grant EP/F003897. The second author was supported in part by Swiss National Science Foundation grant number PP00P2\textunderscore 128557. The third author is supported in part by NSF grant DMS-0856143 and NSF RTG grant DMS-0602191.}
\begin{abstract}
We show that the asymptotic growth rate for the minimal cardinality of a set of simple closed curves on a closed surface of genus $g$ which fill and pairwise intersect at most $K\ge 1$ times is $2\sqrt{g}/\sqrt{K}$ as $g \to \infty$ . We then bound from below the cardinality of a filling set of systoles by $g/\log(g)$. This illustrates that the topological condition that a set of curves pairwise intersect at most once is quite far from the geometric condition that such a set of curves can arise as systoles.
\end{abstract}
\maketitle


\section{Introduction}

We study the size of a collection of simple closed curves on a closed, orientable surface of genus $g\ge 2$ which fill the surface and pairwise intersect at most $K$ times for some fixed $K \geq 1$. In particular, we give upper and lower bounds on the minimum numbers of curves in such a collection. 

The original motivation for this work arose when considering fillings sets of {\it systoles} on closed hyperbolic surfaces (definitions are given in Section \ref{sec:def&note}).  Thurston \cite{Thurston} suggested that the subspace $X_g$ of Teichm\"uller space consisting of those surfaces that admit a filling set of systoles might form a mapping class group invariant spine. However still very little is known about this set of surfaces, and whether $X_g$ is indeed contractible remains to be shown. Another natural but difficult problem is to determine the smallest size of a filling set of systoles. A first approach is to consider the topological constraints on such a set, for it is straightfoward to observe that the curves in a filling set of systoles are nonseparating and can only pairwise intersect at most once. 

In our main result, we are thus generalizing these topological restrictions on filling sets of sytoles by allowing a larger uniform bound on pairwise intersection for filling sets of simple closed curves.

\begin{thm}\label{smallsets}
Let $S$ be a closed, orientable surface of genus $g\ge 2$. The number $n$ of curves in a filling set of simple closed  curves which pairwise intersect at most $K \geq 1$ times satisfies

\begin{equation}\label{inequality}
n^2 - n \geq \frac{4g - 2}{K}
\end{equation}
Moreover, if $N$ is the smallest integer satisfying this inequality, then there exists a set of no more than $N +1$ such curves on $S$.
\end{thm}

Note that Theorem \ref{smallsets} is only really interesting where $g$ is much larger than $K$. We have the resulting corollary: 

\begin{cor}\label{asymptotics}
Let $S$ be a closed, orientable surface of genus $g\ge 2$. The number $n$ of curves in a smallest filling set of simple closed curves which pairwise intersect at most $K\ge 1$ times satisfies $n \sim \frac{2\sqrt{g}}{\sqrt{K}}$; in other words, 
\[ \lim_{g \to \infty} \frac{n}{\sqrt{g}} = \frac{2}{\sqrt{K}} \]
\end{cor}

The topological approximation given by Theorem \ref{smallsets} for the smallest number in a filling set of systoles turns out to be very crude. Our second result says that the number of systoles in a smallest filling set grows with $g$ at a rate of order strictly greater than $\sqrt{g}$.

\

\begin{thm}\label{prop:sys} Let $S$ be a closed, orientable hyperbolic surface of genus $g\ge 2$ with a filling set of systoles $\{ \sigma_1,\hdots,\sigma_n\}$. Then 
$$
n\geq \pi \frac{\sqrt{g(g-1)}}{\log (4g-2)}.
$$
Furthermore there exist hyperbolic surfaces of genus $g$ with filling sets of $n\leq 2g$ systoles.
\end{thm}

Our paper is organized as follows. We begin with definitions and notation in Section \ref{sec:def&note}. For the proof in Section \ref{sec:thm1} of Theorem \ref{smallsets}, we use an Euler characteristic argument to obtain a lower bound on the number in a filling set of curves which pairwise intersect at most $K$ times. To find an upper bound, we give a construction of such a set of curves whose cardinality is at most one larger than the lower bound. As it must work for all $g$ and satisfy this small cardinality condition, this construction is rather cumbersome to describe. Thus in Section \ref{hairy} we produce for $K=1$ easier examples of such small sets of curves whose growth rate, though larger, still has order $\sqrt{g}$. We give an explicit proof that these sets cannot be realized as systoles in Proposition \ref{prop:notsys}. Finally, in Section \ref{sec:growth} we prove Theorem \ref{prop:sys}.

\

 \noindent {\it Acknowledgments.} The authors are thankful to Greg McShane for several helpful discussions which led to the use of the hairy torus example in Section 3. The authors are also very grateful for the anonymous referee's careful reading and useful suggestions. The third author is thankful to Stanford University and the School of Mathematics at the University of Southampton for their hospitality during the conception of this paper.

\section{Preliminaries}\label{sec:def&note}

The results in this paper are stated only for closed, orientable surfaces, although it is possible to extend them to include punctured surfaces. On the other hand, our proofs work explicitly and extensively with compact surfaces with nonempty boundary, and so we begin with some conventions for this purpose.

Let $F$ be a compact, orientable surface with nonempty boundary.  An ${\it arc}$ on $F$ is a simple arc whose endpoints lie on the boundary of $F$ and whose interior lies in the interior of $F$ (that is, we assume that all arcs are properly embedded). We usually confuse an arc with its homotopy class, relative to its basepoints. 

Let $A$ be an arc on the surface $F$. We record the three possibilities for the closure $\overline{F\setminus A}$ of the complement of $A$ in $F$:
\begin{enumerate}
\item The endpoints of $A$ lie in distinct boundary components of $F$. In this case, $A$ is said to be {\it boundary cutting:} the connected surface $\overline{F\setminus A}$ has the same genus as $F$, but one fewer boundary components.
\item The endpoints of $A$ lie in the same boundary component of $F$, and $\overline{F\setminus A}$ is connected. In this case, $A$ is said to be {\it genus cutting:} the surface $\overline{F\setminus A}$ has genus one less than the genus of $F$ and one more boundary component.
\item The endpoints of $A$ lie in the same boundary component of $F$, and $\overline{F\setminus A}$ is disconnected. In this case, $A$ is {\em separating;} the surface $\overline{F\setminus A}$ has genus equal to the genus of $F$ and one more boundary component.
\end{enumerate}
These three possibilities will be revisited in the proof of Theorem \ref{smallsets}.  Note that the removal of an arc of type (1) or (2) does not disconnect the surface; in these two cases, if $F$ is connected, then $\overline{F\setminus A}$ is also connected. The arcs of type (3) occur only towards the end of the proof of Theorem  \ref{smallsets}, when describing the ``last curves'' in a filling set. Note that parts of the boundary of $\overline{F\setminus A}$ come equipped with identifications corresponding to $A$. 

Let $S$ be a closed, orientable surface of genus $g \geq 2$.  A {\em curve} on $S$ is a simple, closed, essential (homotopically nontrivial) curve. We often confuse a curve with its homotopy class; when considering curves which intersect, we generally aim to choose representatives for their homotopy classes with minimal intersection number. 

Recall that a set of curves on $S$ {\it fills} if the complement of the union of its curves consists of simply connected components. Given $K \geq 1$, a set of curves forms a {\it$K$-filling set} if it fills the surface, and if any two curves in the set intersect at most $K$ times. 

By a {\it hyperbolic surface}, we mean a surface equipped with a hyperbolic metric. Every homotopy class of curves on a hyperbolic surface contains a unique geodesic representative. A {\em systole} on a hyperbolic surface is a shortest simple closed geodesic on the surface. A straightforward cut-and-paste argument shows that a filling set of systoles on a surface is a $1$-filling set, and further that all curves in a filling set of systoles are nonseparating.


\section{Proof of Theorem \ref{smallsets}}\label{sec:thm1}

As the title suggests, in this section we give a proof of our main theorem, Theorem \ref{smallsets}.  For this we require lower and upper bounds for the size of a $K$-filling set of simple close curves. 


\subsection{The lower bound}
Let $\gamma_1, \ldots, \gamma_n$ be a set of $n$ curves satisfying the hypotheses of Theorem \ref{smallsets} for some fixed $K \geq 1$. Let $S(m)$ denote the compact surface whose interior is the complement of the union of the first $m$ curves $\gamma_1, \ldots, \gamma_m$ on $S$. Then the Euler characteristic $\CX(S(n))$ of $S(n)$ is at least one.

On $S(m)$, the curve $\gamma_{m+1}$ consists of the union of at most $Km$ arcs whose endpoints correspond to intersections of $\gamma_{m+1}$ with the first $m$ curves.  Deleting a single arc increases the Euler characteristic by at most one. Removing all arcs of $\gamma_{m+1}$, the Euler characteristic of the resulting surface $S(m+1)$ increases by at most $Km$. Thus we have: 
\begin{eqnarray*}
1 \leq \CX(S(n)) & \leq & \CX(S(n-1)) + K(n-1) \\
			& \leq & \CX(S(n-2)) + K(n-2 + n-1) \\
			& \vdots &	\\
			& \leq & \CX(S(1)) + K(1 + \cdots + n-2 + n-1) \\
			& = & 2-2g + K\frac{n(n-1)}{2}		
\end{eqnarray*}
Rearranging the terms gives the inequality of Theorem \ref{smallsets}.


\subsection{The upper bound}
For this we describe a  $K$-filling set with cardinality $N$ or $N+1$, described  inductively in that the $(m+1)^{{\rm st}}$ curve $\gamma_{m+1}$ is given in terms of the previous $\gamma_1,\ldots, \gamma_m$.  As should be apparent in what follows, a minimal such set is in no way unique, and the number of such sets (up to homeomorphism) grows with genus. Part of the difficulty here is in making a choice for these curves in such a way that they can be described systematically, and it should be noted that we make no attempt to keep track of all possible sets of small cardinality.

The curves described will in fact make up an $M$-filling set, where $M= \min\{K,g\}$. This takes into account the relative sizes of $K$ and $g$; for in the case that $g$ is small relative to $K$, the set is $g$-filling so that pairwise intersections between curves are fewer than required. This resolves difficulties in the construction for these cases.  \\


\subsubsection{Notation}
Suppose that $\gamma_1, \ldots, \gamma_m$ are the first $m$ curves chosen for our collection. Recall that $S(m)$ is the compact surface with boundary having interior $S \bs ( \gamma_1\cup \cdots\cup \gamma_m)$. The natural map 
\[ f_m: S(m) \lra S \]
is the identity when restricted to the interior of $S(m)$ and identifies boundary components corresponding to the curves $\gamma_1, \ldots, \gamma_m$. For most cases (see Section \ref{intermediatearcs} below), the curve $\gamma_{m}$ is described by $M(m-1)$ arcs $\{ A_{m}^{i,j} \}_{1 \leq j \leq M}^{1 \leq i \leq m-1}$ on $S(m-1)$ such that 
\[ \gamma_{m} = f_{m-1}\bigg(\bigcup_{i=1}^{m-1} \bigcup_{j=1}^M A_{m}^{i,j}\bigg). \]
If $\gamma_m$ is one of the ``last curves'', then it is given by a similar union, but of possibly fewer than $M(m-1)$ arcs (see Section \ref{lastcurves} below). Denote by $b_m^{i,j}$ and $t_m^{i,j}$ the base and terminal points of the arc $A_m^{i,j}$, respectively. The order in which the arcs are joined will be clear when discussing the relationship between the base and terminal points of arcs.  For brevity and clarity, we say that the arc $A_m^{i,j}$ is {\it genus cutting} (respectively, {\it boundary cutting}) if it is genus cutting (respectively, boundary cutting) in the surface whose interior is:
\[ S(m-1) \ \backslash \  \bigcup_{k=1}^{i-1} \bigcup_{\ell=1}^{M} A_m^{k,\ell} \ \backslash \   \bigcup_{\ell=1}^{j-1} A_m^{i,\ell}  \]

\


\subsubsection{The first two curves}
The first curve $\gamma_1$ is a nonseparating simple closed (oriented) curve on $S$. To describe the curves $\gamma_{m}$ for $m \geq 2$, we first consider when the genus of $S(m-1)$ is $M$ more than the number of boundary cutting arcs from $\gamma_{m-1}$. If this is already not the case for $m=2$, the remaining curves on $S$ are ``last curves,'' as described in Section \ref{lastcurves} below.

Let $\gamma_1^+$ (respectively, $\gamma_1^-$) be the component in $S(1)$ of $f_1^{-1}(\gamma_1)$ for which the surface lies on the left (respectively, right). Choose $M$ distinct points $p_1, p_2, \ldots, p_M$ on $\gamma_1$,  ordered for convenience according to the orientation of $\gamma_1$. For $1 \leq k \leq M$, let $p_k^+$ be a point on $\gamma_1^+$, and $p_k^-$ a point on $\gamma_k^-$ such that $f_1(p_k^+) = f_1(p_k^-) = p_k$.

For $1 \leq j < M$, the base point of $A_2^{1,j}$ is $b_2^{1,j} = p_j^+$, and the terminal point is $t_2^{1,j} = p_{j+1}^-$; the terminal point $t_2^{1,M}$ of $A_2^{1,M}$ is $p_1^-$. The union $\bigcup_{j=1}^{M} A_2^{1,j}$ then maps under $f_1$ to a simple closed curve $\gamma_2$.  We choose the arcs so that $A_2^{1,1}$ is boundary cutting and each remaining arc $A_2^{1,j}$, $2 \leq j \leq M$, is genus cutting.  Hence, the surface $S(2)$ will have $M$ boundary components.

\


\subsubsection{The intermediate curves}\label{intermediatearcs}
The curve $\gamma_{m+1}$, $m \geq 2$, is considered an ``intermediate curve'' as long as the genus of $S(m)$ is at least $M$ greater than the number of arcs $A_m^{i,j}$ which are boundary cutting in $S(m-1)$. If $S(m)$ has smaller genus, then $\gamma_{m+1}$ is described as one of the ``last curves'' below.

For each arc $A_m^{i,j}$ of $\gamma_m$, there will be a corresponding arc $A_{m+1}^{i,j}$ of $\gamma_{m+1}$. If $A_m^{i,j}$ is genus cutting, then $A_{m+1}^{i,j}$ is chosen to be boundary cutting; likewise if $A_m^{i,j}$ is boundary cutting, then $A_{m+1}^{i,j}$ is genus cutting. This convention keeps down the number of boundary components and connected components of the surface $S(m)$.

\

We now define the base and terminal points of the arcs $A^{i,j}_{m+1}$ on $S(m)$.  The first arc $A_{m+1}^{1,1}$ is based at a point $b_{m+1}^{1,1}$ in the interior of a component of $f_m^{-1}(\gamma_1)$ which, when projected to $S(m-1)$, contains the base point $b_m^{1,1}$ of $A_m^{1,1}$.   Once we define the terminal points, all of the base points are determined because the union $\bigcup_{i=1}^{m} \bigcup_{j=1}^M A_{m+1}^{i,j}$ mapping under $f_m$ to the closed curve $\gamma_{m+1}$
forces relations between base points and terminal points of arcs. Namely, when $i=1$ and $1 < j \leq M$, we define the arcs $A^{1,j}_{m+1}$ to be based at the points 
\[ b_{m+1}^{1,j} = f_m^{-1}(f_m(t_{m+1}^{1,j-1})) \bs \{ t_{m+1}^{1,j-1} \} \]

\noindent Then for $1 < i \leq m$ and $j=1$, we define the arcs $A_{m+1}^{i,1}$ to be based at the points 
\[ b_{m+1}^{i,1} = f_m^{-1}(f_m(t_{m+1}^{i-1,M})) \bs \{ t_{m+1}^{i-1,M} \} \]
Finally, we define the arcs $A_{m+1}^{i,j}$, for $1 < i \leq m$ and $1 < j \leq M$, to be based at the points 
\[ b_{m+1}^{i,j} = f_m^{-1}(f_m(t_{m+1}^{i,j-1})) \bs \{ t_{m+1}^{i,j-1} \} \]
We now describe the terminal points for these arcs; these will depend on whether the arc in question is genus cutting or boundary cutting. 

\

We first describe how to terminate the arcs $A_{m+1}^{i,j}$ for $1\leq i \leq m-1$, $1 \leq j \leq M$, and $(i,j) \neq (m-1,M)$. The notation $\mathcal{N}(\gamma_m)$ will always indicate a thin collar neighbourhood of $\gamma_m$ which does not contain the intersection of any two curves from the set $\{ \gamma_1, \ldots, \gamma_{m-1}\}$. For some such neighbourhood the arc $A_{m+1}^{i,j}$ terminates at a point $t_{m+1}^{i,j}$ contained in the same component of $f_m^{-1}(\mathcal{N}(\gamma_m))$ as one of the two points of $f_m^{-1}(f_{m-1}(t_m^{i,j}))$. The terminal point $t_{m+1}^{i,j}$ lies on a component of $f_m^{-1}(\gamma_k)$, where $k=i+1$ if $M=1$, where $k = i$ if $1 \leq j < M$, and where $k=i+1$ if $j=M$.

Suppose that $A_m^{i,j}$ is genus cutting so that $A_{m+1}^{i,j}$ is boundary cutting. Then $A_{m+1}^{i,j}$ terminates at a point $t_{m+1}^{i,j}$ lying in a different boundary component of $S(m)$ from the base point $b_{m+1}^{i,j}$. On the other hand, if $A_m^{i,j}$ is boundary cutting so that $A_{m+1}^{i,j}$ is genus cutting, then the terminal point $t_{m+1}^{i,j}$ of $A_{m+1}^{i,j}$ is on the same boundary component of $S(m)$ as the base point $b_{m+1}^{i,j}$.

\

This leaves to describe those arcs of $\gamma_{m+1}$ which intersect the arcs of the curve $\gamma_m$. Let $q_1,\ldots,q_M$ in $S(m-1)$ be $M$ points along the arc $A_m^{m-1,M}$ ordered according to the orientation of $A_m^{m-1,M}$. Again $\mathcal{N}(\gamma_m)$  can be chosen so that there is a unique point $q_i^+ \in S(m)$ of $f_m^{-1} (f_{m-1}(q_i))$ which lies in the same component of $f_m^{-1}(\mathcal{N}(\gamma_m))$ as 
\[ f_m^{-1}(f_m (t_{m+1}^{m-1,M-1})) \bs \{ t_{m+1}^{m-1,M-1} \} \]
and $q_i^-$ the other point in $f_{m}(f_{m-1}^{-1}(q_i))$. The arc $A_{m+1}^{m-1,M}$ terminates at $t_{m+1}^{m-1,M} = q_1^-$, and each arc $A_{m+1}^{m,j}$ for $1 \leq j \leq M-1$ terminates at $t_{m+1}^{m,j} = q_{j+1}^-$. The arcs $A_{m+1}^{m,j}$, $1 \leq j \leq M$, are all genus cutting, with the final arc $A_{m+1}^{m,M}$ terminating at 
\[ t_{m+1}^{m,M} = f_m^{-1}( f_m(b_{m+1}^{1,1})) \bs \{ b_{m+1}^{1,1} \} \]
closing up the curve $\gamma_{m+1}$. 

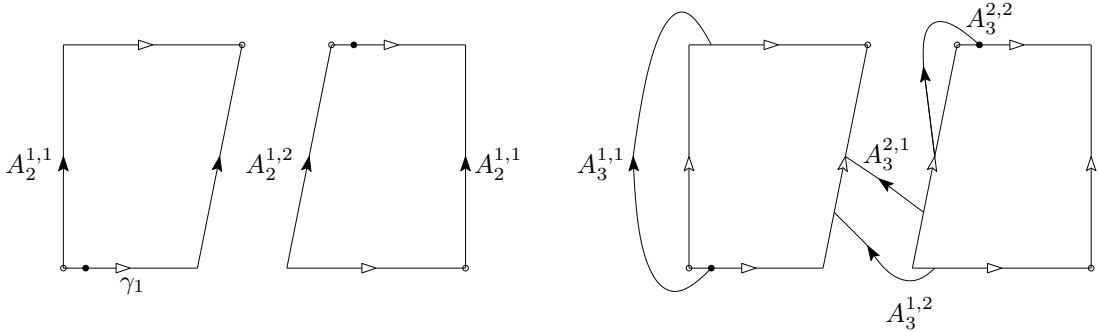
\begin{figure}[h]
\leavevmode \SetLabels
\L (.1*.1) $\gamma_{1}$\\
\L (0*.45) $A_2^{1,1}$\\
\L (.41*.45) $A_2^{1,1}$\\
\L (.21*.45) $A_2^{1,2}$\\
\L (.5*.45) $A_3^{1,1}$\\
\L (.77*0) $A_3^{1,2}$\\
\L (.75*.48) $A_3^{2,1}$\\
\L (.84*.89) $A_3^{2,2}$\\
\endSetLabels
\par
\begin{center}
\AffixLabels{\centerline{\epsfig{file =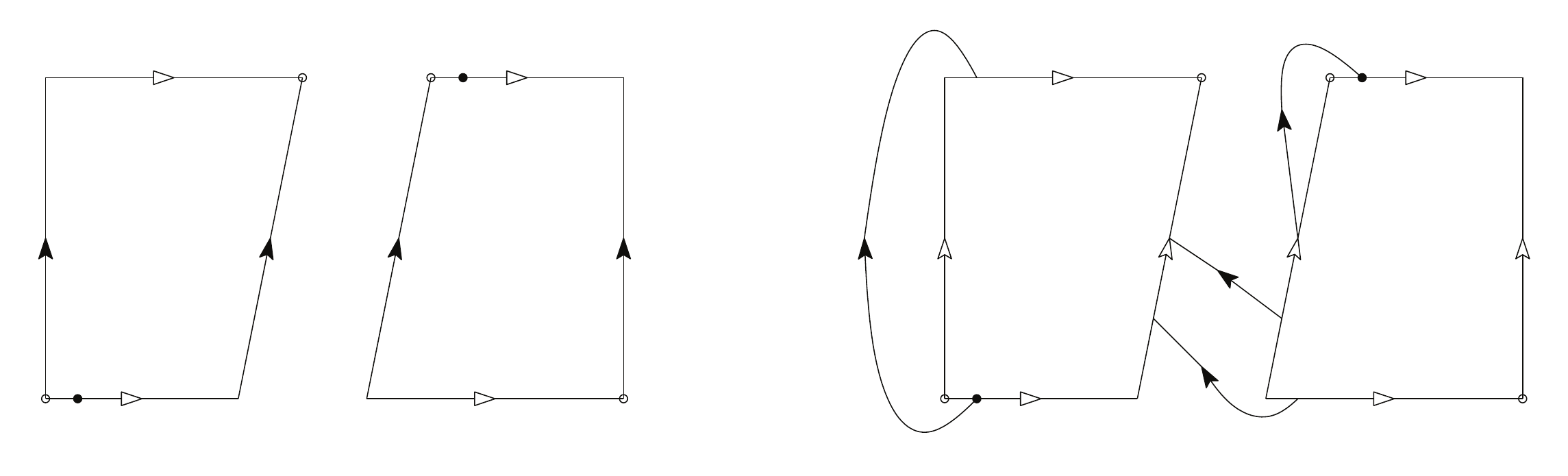, width=14.5cm, angle= 0}}}
\end{center}
\caption{On the left is the boundary of $S(2)$ for $K=2$ and a surface of genus $g\geq 5$. On the right, the arcs of the curve $\gamma_3$ are shown on $S(2)$. The open points indicate the basepoint of $\gamma_2$, while the solid points indicate the basepoint of $\gamma_3$. The arcs $A_2^{1,1}$ and $A_3^{1,2}$ are boundary cutting, while the arcs $A_2^{1,2}$, $A_3^{1,1}$, $A_3^{2,1}$, and $A_3^{2,2}$ are genus cutting.}
\label{fig:example}
\end{figure}

\


\subsubsection{The last curves}\label{lastcurves}
Finally we consider the case when the genus of $S(m)$, $m \geq 2$, is not large enough to define $Mm$ nonseparating arcs as described above. This occurs when the genus of $S(m)$ is less than the number of boundary cutting arcs $A_m^{i,j}$ in $S(m-1)$, or else exceeds it by less than $M$. We treat these two cases separately. 

\

\noindent Case (1): The genus of $S(m)$ is at most the number of boundary cutting arcs $A_m^{i,j}$. We choose arcs $A_{m+1}^{1,1}, A_{m+1}^{1,2}, \ldots, A_{m+1}^{1,M}, A_{m+1}^{2,1}, A_{m+1}^{2,2}, \ldots$ as above until we have that complement of 
\[ \left( \bigcup_{1 \leq k < i_1} \bigcup_{1 \leq j \leq M} A_{m+1}^{k,j} \right)
\cup \left( \bigcup_{1 \leq \ell < j_1}A_{m+1}^{i_1,\ell} \right) \] 
in $S(m)$ is a sphere with holes, and so that $A_m^{i_1,j_1}$ is a boundary cutting arc of $\gamma_m$. For $A_m^{i_1,j_1}$ and every subsequent boundary cutting arc $A_m^{i,j}$ of $\gamma_m$, the corresponding arc $A_{m+1}^{i,j}$ of $\gamma_{m+1}$ is chosen to lie entirely in a component of some thin collar neighbourhood $\mathcal{N}(\gamma_m)$ of $\gamma_m$; the arc $A_{m+1}^{i,j}$ is then separating, with the property that one of the components of its complement is simply connected.

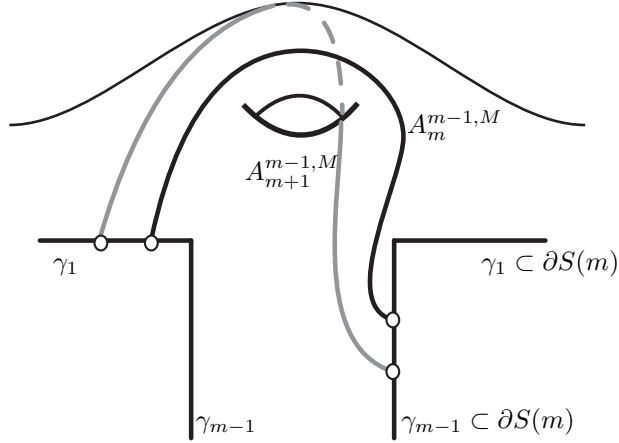
\begin{figure}[h]
\leavevmode \SetLabels
\L (.458*.60) $A_{m+1}^{m-1,M}$\\
\L (.6*.69) $A_m^{m-1,M}$\\
\L (.415*.13) $\gamma_{m-1}$\\
\L (.595*.13) $\gamma_{m-1}\subset \partial S(m)$\\
\L (.665*.43) $\gamma_1\subset \partial S(m)$\\
\L (.29*.43) $\gamma_1$\\
\endSetLabels
\par
\begin{center}
\AffixLabels{\centerline{\epsfig{file =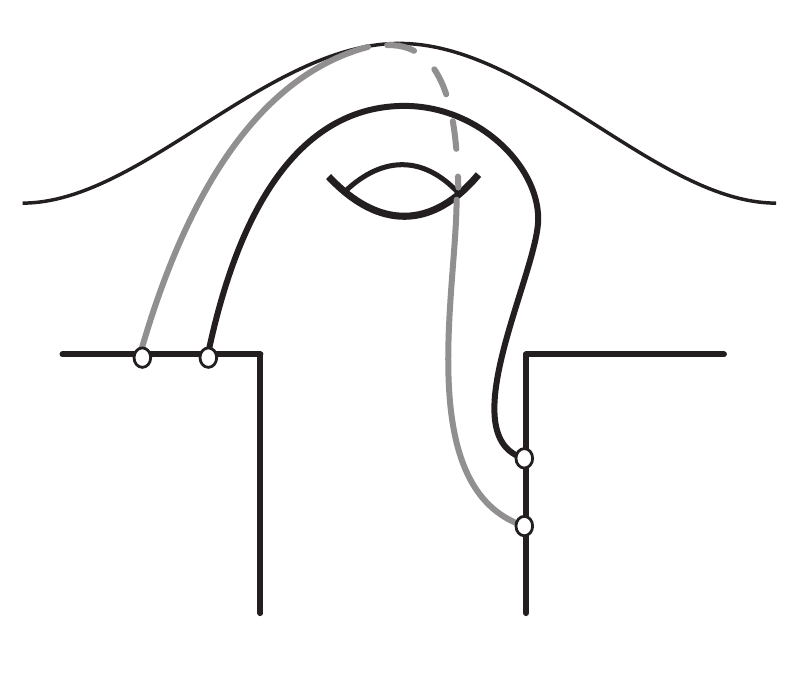, width=8cm, angle= 0}}}
\end{center}
\caption{Here $\gamma_{m+1}$ closes up with the arc $A_{m+1}^{m-1,M}.$}
\label{fig:case1}
\end{figure}

\begin{figure}[h]
\leavevmode \SetLabels
\L (.545*.51) $A_{m+1}^{m,1}$\\
\L (.46*.60) $A_{m+1}^{m-1,M}$\\
\L (.415*.3) $A_m^{m-1,M}$\\
\L (.557*.13) $\gamma_{m-1}\subset \partial S(m)$\\
\L (.665*.42) $\gamma_1\subset \partial S(m)$\\
\endSetLabels
\par
\begin{center}
\AffixLabels{\centerline{\epsfig{file =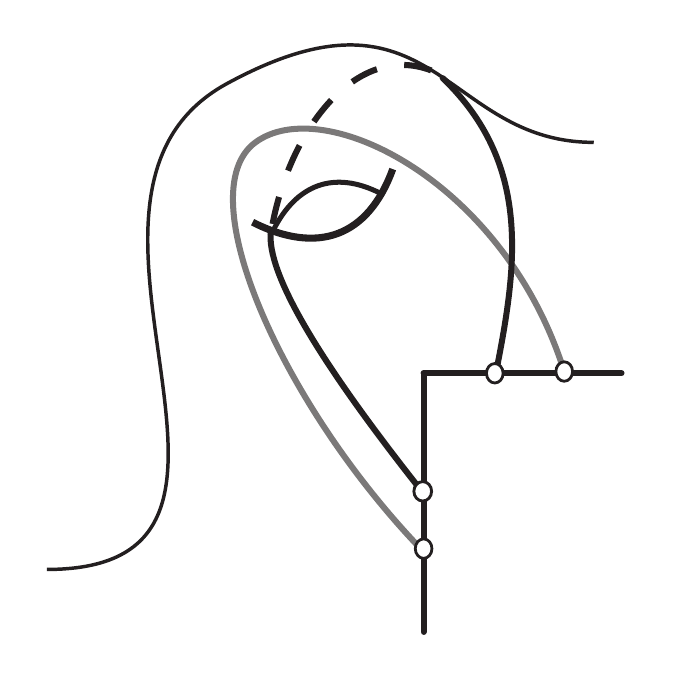, width=7cm, angle= 0}}}
\end{center}
\caption{Here $\gamma_{m+1}$ closes up with the arc $A_{m+1}^{m,1}$.}
\label{fig:case2}
\end{figure}

Now recall that the arc $A_m^{m-1,M}$ is genus cutting; the last arc of $\gamma_{m+1}$ is then a boundary cutting arc we label $A_{m+1}^{m-1,M}$ (not to be confused with $A_{m+1}^{m-1,M}$ as chosen above, where $\gamma_{m+1}$ is an intermediate curve), if such an arc exists (see Figure \ref{fig:case1}). Otherwise $A_{m+1}^{m-1,M}$ is a boundary cutting arc whose terminal point $t_{m+1}^{m-1,M}$ is mapped by $f_m$ to a point on $f_{m-1}(A_m^{m-1,M})$, and a last separating arc $A_{m+1}^{m,1}$ closes up $\gamma_{m+1}$ (see Figure \ref{fig:case2}). (Again, $A_{m+1}^{m,1}$ is chosen according to whether $\gamma_{m+1}$ is an intermediate curve, as in the last section, or a last curve.)

\

\noindent Case (2): The genus of $S(m)$ exceeds the number of boundary cutting arcs $A_m^{i,j}$ by $L < M$. In this case, the genus cutting arcs $A_{m+1}^{m,j}$ for $1 \leq j \leq L-1$ are as described for the intermediate curves. The very last arc of $\gamma_{m+1}$ is the genus cutting arc $A_{m+1}^{m,L}$ which closes up $\gamma_{m+1}$.

\

In either Case (1) or Case (2), the resulting complementary surface (when the arcs described so far have been deleted) has at most one non simply connected component (in Case (2) the surface $S(m+1)$ is connected). If all components of the surface $S(m+1)$ are simply connected, then $\gamma_1, \ldots, \gamma_{m+1}$ is a filling set. Otherwise, the final curve of our filling set is $\gamma_{m+2}$, consisting of the final set of arcs we describe as follows.

\

We begin by noting that the surface $S(m+1)$ has genus $0$ and multiple boundary components come from genus cutting arcs in the previous step. As before we choose arcs $A_{m+2}^{i,j}$ to be boundary cutting if $A_{m+1}^{i,j}$ was genus cutting. Whenever $A_{m+1}^{i,j}$ is either boundary cutting or separating, $A_{m+2}^{i,j}$ is chosen to be a separating in $S(m+1)$. Then if $\gamma_{m+1}$ arose from Case (1) above, the final arc of $\gamma_{m+2}$ is separating. If on the other hand $\gamma_{m+1}$ arose from Case (2) above, then the final arc chosen for $\gamma_{m+2}$ is a boundary cutting arc closing up $\gamma_{m+2}$, if such an arc exists; otherwise a boundary cutting arc terminating on that final arc of $\gamma_{m+1}$ is chosen, followed by a separating arc which closes up $\gamma_{m+2}$. The arcs chosen reduce the number of boundary components or cut off simply connected pieces. The full set of curves is now necessarily filling.

\

Now that we have described the set of curves, we must establish that the curves are nonseparating and pairwise nonhomotopic. Notice that up until at least $k=m$, the curves $\gamma_k$ are nonseparating as the complement of their union has one component. They are also pairwise nonhomotopic as a consequence of the following well known lemma (see, for example, the Bigon Criterion in \cite{Farb-Margalit}).

\begin{lem}\label{lem:nonseparating}
If $\gamma$ and $\delta$ are two simple closed curves on $S$ such that $S\setminus \{\gamma\cup \delta\}$ is connected, then $\gamma$ and $\delta$ intersect minimally (among all representatives of their respective isotopy classes).
\end{lem}

\noindent This takes care of the ``intermediate curves.'' 

The ``last curves,'' $\gamma_{m+i}$, $i=1,2$, are slightly more problematic as they may be composed of arcs some of which separate the surface they are defined on. However by construction these separating arcs can be isotoped (relative endpoints) to arcs on the boundary of their respective surface so that the union of the resulting set of arcs do not separate. These isotopies induce an isotopy of the corresponding curve $\gamma_{m+i}$ on $S$ whose union with the previous curves is not separating; thus the isotoped curve is itself not separating.

We claim that $\gamma_{m+1}$ and $\gamma_{m+2}$ are not homotopic to each other, nor to any of the other curves in the set. The curve $\gamma_{m+i}$, $i=1,2$, necessarily contains the image under $f_{m+i-1}$ of a boundary cutting arc $A$ with endpoints on two boundary components of $S(m+i-1)$. We begin by observing that if $S(m+i-1)$ has more than one boundary component, then they are simple (essential) closed curves on the surface $S$. This is because multiple boundary components comes from genus cutting arcs at the previous step. Thus there is always an arc on a boundary component which on $S$ is glued to an arc of some other boundary component. Consider a simple curve obtained by gluing a simple path between the midpoints of the two arcs: it intersects the boundary component exactly once and by the bigon criterion, both are essential. Now let $\eta$ denote one of them. After a possible isotope of $\eta$, Lemma~\ref{lem:nonseparating} implies that the curve $\gamma_{m+i}$ intersects $\eta$ nontrivially, while the previous curves $\gamma_1, \ldots, \gamma_{m+i-1}$ have empty intersection with $\eta$. This proves the claim.

\subsubsection{The final step} Now recall that $N$ is the smallest integer satisfying the inequality (\ref{inequality}). Note that up through stage $m$, when genus runs out, each arc increases the Euler characteristic of the complementary surface by one. As a consequence, we must have $m+1 \leq N$ (see the argument for the lower bound). On the other hand, only at most two additional curves $\gamma_{m+1}$ and $\gamma_{m+2}$ are required to fill the surface, and so $N \leq m+2$. Therefore $m$ is equal to $N-2$ or $N-1$, and this completes the proof of Theorem \ref{smallsets}.


\section{A simple small filling set}\label{hairy}

The $1$-filling set of curves described in Theorem \ref{smallsets} is not easy to describe or visualize, and so here we provide an alternative family of examples of $1$-filling sets on closed, orientable surfaces. They are not optimally small in size, but still have growth rate of order $\sqrt{g}$. Although purely topological, these examples are inspired by Buser's hairy torus example (see \cite[Section 5.3]{Buser}), used there to find a lower bound on the Bers constant. We shall see quite explicitly in Proposition \ref{prop:notsys} that the family  provides examples of $1$-filling sets which cannot be realized as geometric systoles. 

\begin{figure}[h]
\leavevmode \SetLabels
\L (.33*1.025) $\alpha_k$\\
\L (.33*-.05) $\alpha_{k+1}$\\
\L (.33*.67) $\gamma_k$\\
\L (.515*.23) $\delta_k$\\
\L (.2*.2) $\beta_k$\\
\L (.775*.2) $\beta_{k+1}$\\
\endSetLabels
\par
\begin{center}
\AffixLabels{\centerline{\epsfig{file =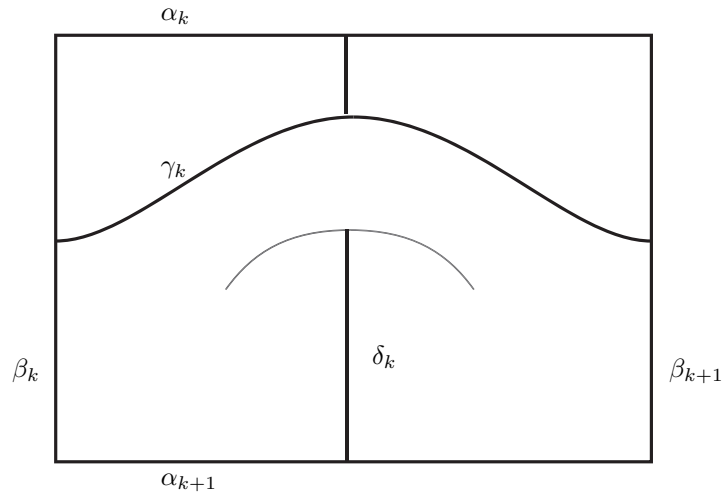, width=8cm, angle= 0}}}
\end{center}
\caption{The torus building block}
\label{fig:square}
\end{figure}

For all integers $N>0$, consider the family of surfaces of genus $N^2+1$. Such surfaces can topologically be constructed by gluing together $N^2$ one holed tori with ``square" boundary (see Figure \ref{fig:square} for the torus with square boundary and Figure \ref{fig:example} to see how they are pasted together).

\begin{figure}[h]
\leavevmode \SetLabels
\L (.35*.92) $\alpha_1$\\
\L (.35*.15) $\alpha_{1}$\\
\L (.26*.79) $\gamma_1$\\
\L (.25*.265) $\gamma_{N}$\\
\L (.32*.74) $\delta_1$\\
\L (.68*.74) $\delta_{N}$\\
\L (.21*.22) $\beta_1$\\
\L (.77*.22) $\beta_{1}$\\
\endSetLabels
\par
\begin{center}
\AffixLabels{\centerline{\epsfig{file =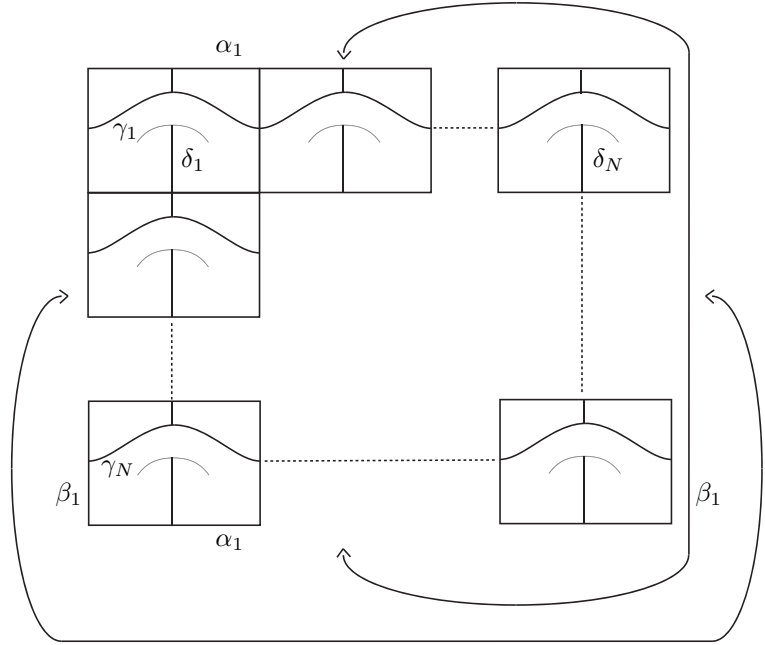, width=10cm, angle= 0}}}
\end{center}
\caption{A surface of genus $N^2+1$}
\label{fig:example}
\end{figure}

Denote the horizontal curves on the surface which were originally the horizontal boundaries of the squares as $\alpha_k$, $k=1,\hdots, N$, and the corresponding vertical ones as $\beta_k$, $k=1,\hdots, N$ (see Figure \ref{fig:square} for labels). Also consider the horizontal and vertical curves $\gamma_k$, $k=1,\hdots, N$ and $\delta_k$, $k=1,\hdots, N$ as in Figures \ref{fig:square} and \ref{fig:example}. Note that the union of all four sets gives a set of $4N$ curves which pairwise intersect at most once and whose complementary region is a set of $N^2$ disks. As the genus is equal to $N^2+1$, the cardinality of the set of curves grows in genus like $4 \sqrt{g}$.\\

The following proposition shows that the curves described above cannot be realized as systoles.

\begin{prop}\label{prop:notsys} For any $N\ge 4$ and any choice of (hyperbolic) metric on a closed, orientable surface of genus $N^2+1$, the curves $\{\alpha_k\}_{k=1}^N$ and $\{\beta_k\}_{k=1}^N$ described above cannot be systoles.
\end{prop}

\begin{proof}
Consider the $N^2$ simple nontrivial curves, say $\epsilon_k$, $k=1,\hdots,N^2$, which were originally the boundary curves of the one holed tori used to construct our surface. Now suppose that there exists a hyperbolic metric for which the curves $\{\alpha_k\}_{k=1}^N$ and $\{\beta_k\}_{k=1}^N$ are systoles of length $x$. Denote the length of a curve $\gamma$ with this metric as $\ell(\gamma)$. Each $\epsilon_k$ is homotopic to a closed curve consisting of two segments of two consecutive $\alpha_{i}$s and two segments of two consecutive $\beta_{j}$s, and one can construct curves homotopic to all of the $\epsilon_k$ using each segment exactly twice. Thus one has the following inequality:

$$
\sum_{k=1}^{N^2} \ell (\epsilon_k) < 2 \sum_{k=1}^{N} (\ell(\alpha_k)+\ell(\beta_k)).
$$
But as $\ell(\epsilon_k) \geq x$ and $\ell(\alpha_k)=\ell(\beta_k)=x$ this gives
$$
N^2 x < 4N x
$$
which is a contradiction for $N \ge 4$.
\end{proof}
We remark that the proposition also holds for $N=2,3$; however these cases require separate consideration from $N \geq 4$. As they serve no real purpose in this paper, their proof is omitted.


\section{Growth of filling systoles}\label{sec:growth}

As we have seen in the previous section, the $K$-filling sets and the sets of curves which are realizable as systoles are distinct. Here we investigate the minimal cardinality of filling sets of systoles.

\begin{named}{Theorem \ref{prop:sys}}
Let $S$ be a closed, orientable hyperbolic surface of genus $g\ge 2$ with a filling set of systoles $\{ \sigma_1,\hdots,\sigma_n\}$. Then 
$$
n\geq \pi \frac{\sqrt{g(g-1)}}{\log (4g-2)}.
$$
Furthermore there exist hyperbolic surfaces of genus $g$ with filling sets of $n\leq 2g$ systoles.
\end{named}

\begin{proof}
We begin with the first inequality. The set of systoles fills, and so their complementary region on the surface is a union of hyperbolic disks with piecewise geodesic boundary. The total perimeter of the set of disks obtained this way is twice the sum of the lengths of the systoles. It is easy to see that the length $x$ of a systole of a surface of genus $g$ satisfies the inequality $x\leq 2 \log( 4g-2)$ (see \cite[Lemma 5.2.1]{Buser}). So the total perimeter of the disks does not exceed $L\leq 4n \log(4g-2)$.\\

Now the total area of the disks equals the area of the surface ($A=2\pi(2g-2)$). By the isoperimetric inequality in the hyperbolic plane we have 

$$
L \geq \sqrt{A^2+4\pi A}
$$
which implies
$$
4n \log(4g-2) \geq 4\pi \sqrt{g(g-1)}
$$
and the desired inequality follows.\\

To show the second inequality, it suffices to exhibit an example of a surface with a filling set of systoles of cardinality $2g$. For $g\geq 2$, one can construct a surface by pasting together, in an appropriate way, four copies of a regular right angled $(2g+2)$-gon. Note that for every $g$, it is straightforward to establish that there is only one such polygon; we denote its (uniquely determined) side length $\ell$, and let $m$ denote the barycenter of the vertices.

By standard hyperbolic trigonometry on the polygon, we obtain that

$$
\ell=2 \arcsinh\sqrt{\cos\frac{\pi}{g+1}}.
$$
Note that the distance in the interior of the polygon between any two nonadjacent edges is greater than $\ell$ (and strictly greater if the two edges are more than one apart from being adjacent). To see this, consider any two such edges $e_1$ and $e_2$, and the two paths, say $c_1$ and $c_2$, of length $d_g$ between $m$ and the midpoint of these edges. The paths $c_1$ and $c_2$ plus the distance realizing path between $e_1$ and $e_2$ are part of a pentagon with right angles, except for the angle $\frac{\pi}{g+1} k$ between $c_1$ and $c_2$, where $k$ is the minimal number of edges between $e_1$ and $e_2$. Note that this pentagon is uniquely determined up to isometry by $g$ and $k$. Finally: for $k=1$ one obtains the distance $\ell$ calculated above, and the distance between $e_1$ and $e_2$ is a strictly increasing function in $k$.\\

In pairs, one glues the polygons along every second edge to obtain two $(g+1)$-holed spheres of boundary lengths all equal to $2\ell=4 \arcsinh\sqrt{\cos(\frac{\pi}{g+1})}$. These spheres have a natural set of ``seams" corresponding to the edges of the polygons that have been pasted together. Now one glues along the boundaries of the two surfaces such that the two end points of each seam are pasted exactly to the two end points of a same seam on the other holed sphere. There are now $2g+2$ filling curves on the resulting surface of length $2\ell$. It suffices to see that these curves are systoles. To see this, note that any simple nontrivial simple closed curve not homotopic to one of the $2g+2$ curves constructed above necessarily contains at least two arcs which join nonadjacent edges on one of the four polygons. As any such arc has length more than $\ell$, this implies that any nontrivial simple closed curve has length more than $2\ell$.\\

The surface we have constructed has $2g+2$ systoles but it is easy to see that any subset of these systoles $\sigma_1,\hdots,\sigma_{2g}$ with intersection number $\ii(\sigma_k,\sigma_{k+1})=1$ for $k=1$ to $2g-1$ fills, and that their complementary region is a (right angled) connected polygon.
\end{proof}

Observe that the lower bound in the previous theorem implies that there are at most a finite number of genera for which the 1-filling sets constructed in the proof of Theorem \ref{smallsets} are realizable as systoles.\\

Let us conclude by remarking the following. Families of surfaces with number of filling systoles of growth order $\frac{g}{\log(g)}$ are possibly quite difficult or impossible to find. The proof of Theorem \ref{prop:sys} shows that if such a family exists then the lengths of the systoles must grow at least like $\log(g)$. Although such families of surfaces are known to exist (the first examples were found in \cite{BuserSarnak}), they are very difficult to construct and most known examples come from using arithmetic groups in some way. Furthermore, the order of growth of the number of systoles for these families is at least linear. One would have to analyse the topology of the sets of systoles of such surfaces to see if there are filling subsets with less than linear growth.



\begin{thebibliography}{99}

\bibitem{Buser}
Peter Buser, {\it Geometry and Spectra of Compact Riemann Surfaces}, Progress in Mathematics, 106 Birkh\"auser Boston, Inc., Boston, MA, 1992.

\bibitem{BuserSarnak}
Peter Buser and Peter Sarnak, {\it On the period matrix of a {R}iemann surface of large genus}, Invent. Math. 117(1), 106, 27--56, 1994.

\bibitem{Farb-Margalit}
Benson Farb and Dan Margalit, {\it A Primer on Mapping Class Groups}. 
\newline {\tt http://www.math.utah.edu/$\sim$margalit/primer/}.

\bibitem{Thurston}
William Thurston, {\it A spine for Teichm\"uller space}, preprint.


\end{thebibliography}
\end{document}